\documentclass [a4paper,english]{article}
\usepackage{amsthm, amsmath, amssymb, amsbsy, amsfonts, latexsym, euscript}
\usepackage[cp1251]{inputenc}
\usepackage[T2A]{fontenc}
\usepackage{indentfirst,epigraph}
\usepackage[english,russian]{babel}

\def\le{\leqslant}
\def\leq{\leqslant}
\def\ge{\geqslant}
\def\geq{\geqslant}
\def\phi{\varphi}

\def\bar{\overline}

\def\kappa{\varkappa}

\newtheorem{theorem}{\bf Теорема}

\newtheorem{corollary}{\bf Следствие}

\theoremstyle{remark}
\newtheorem{remark}{\bf Замечание}
\newtheorem{definition}{\bf Определение}[section]

\title{Inequalities  of M.G. Krein, Yu.V. Linnik and E.A. Gorin
 for positive definite functions}
\author{A.B. Pevnyi ,  S.M. Sitnik}

\date{}

\begin{document}
\maketitle

\selectlanguage{english}
\begin{abstract}

We investigate inequalities of M.G. Krein, Yu.V. Linnik and E.A. Gorin for positive definite functions. Modifications and generalizations of these inequalities are proved. We also prove that multipoint E.A. Gorin's inequality follows from two--point  M.G. Krein's inequality.

\end{abstract}

\newpage

\selectlanguage{russian}

\begin{center}
{\Large \bf
НЕРАВЕНСТВА М.~Г.~КРЕЙНА, Ю.~В.~ЛИННИКА И Е.~А.~ГОРИНА ДЛЯ ПОЛОЖИТЕЛЬНО ОПРЕДЕЛЁННЫХ ФУНКЦИЙ}
\end{center}

\begin{abstract}
Исследуются неравенства М.~Г.~Крейна, Ю.~В.~Линника и Е.~А.~Горина для положительно определённых функций, заданных на вещественной оси. Доказываются модификации и обобщения этих неравенств. Показано, что многоточечное неравенство Е.~А.~Горина можно вывести из двухточечного неравенства М.~Г.~Крейна.
\end{abstract}

\vskip10mm

\noindent

\section{Введение}

Теория положительно определённых функций (п.о.ф.) возникла в начале 20 века
на стыке нескольких разделов математики.
Из литературы по п.о.ф. отметим одну из первых оригинальных
работ \cite{Mat}, содержащую по существу все основные современные определения, монографию Бхатиа \cite{Bha}, обзор Стюарта \cite{Ste}.

В данной работе  доказываются варианты известных неравенств М.Г.
Крейна, Ю.В. Линника и Е.А. Горина. Доказаны некоторые обобщения этих неравенств, установлена их взаимосвязь. Отметим, что для авторов инициирующей послужила статья Е.А. Горина \cite{Gor1}.

Дадим основное определение.

\begin{definition} Функция $f:\mathbb{R}\to \mathbb{C}$
называется \textit{положительно определённой функцией} (п.о.ф.), если для любого N,
любых $x_{1},\mathellipsis x_{N}\in \mathbb{R}$ и любой последовательности
\textit{комплексных} чисел $z_{1},\mathellipsis ,z_{N}\in \mathbb{C}$ выполнено неравенство
\begin{equation}
\label{eq1}
\sum\limits_{k,j=1}^N {f\left( x_{k}-x_{j} \right)z_{k}\bar{z_{j}}\ge 0.}
\end{equation}
\end{definition}

Из \eqref{eq1} следует, что $|f(x)| \leq f(0)$ и $f(-x)=\bar{f(x)}$ для всех $x\in\mathbb{R}$. Вещественная часть п.о.ф. снова является п.о.ф. По поводу этих фактов см. \cite{Bha}--\cite{Ste}, \cite{ChLi}, \cite{Zas}. Важными примерами п.о.ф. являются $\exp(ix), \cos(x),  \exp(-x^2)$.

\section{Неравенство М.Г. Крейна и его обобщение}

Справедливо классическое неравенство  М.Г. Крейна \cite{Kre}:
\begin{equation}
\label{eq2}
{[f\left( x \right)-f(y)]}^{2}\le 2f\left( 0 \right)\Re \left[ f\left( 0
\right)-f\left( x-y \right) \right];\, x,y\in \mathbb{R}.
\end{equation}

Установим теперь новый результат --- обобщённое неравенство М.Г. Крейна, содержащее произвольный комплексный параметр $\alpha$.

\begin{theorem} Пусть функция $f(x)$ есть п.о.ф. на $\mathbb{R}$. Тогда для любых $x,y\in \mathbb{R}$ и любого $\alpha\in\mathbb{C}, |\alpha|=1$ справедливо неравенство
\begin{equation} \label{eq3}
|\alpha f(x) - f(y)|^2 \leq 2 f(0) \Re\left[f(0)-{\alpha f(x-y)}\right].
\end{equation}
\end{theorem}

Доказательство. Не умаляя общности, можно считать, что $f(0)=1$. Выберем в \eqref{eq1} три точки $0,x,y$ и рассмотрим матрицу

\[A=
\left( {\begin{array}{*{20}c}
	1 & \bar{f(x)} & \bar{f(y)}\\
	f(x) & 1 & f(x-y)\\
	f(y) & f\bar{(x-y)} & 1\\
	\end{array} } \right).
\]
По условию \eqref{eq1} $(Az,z)\geq 0$ для любого $z\in\mathbb{C}^3$.
Выберем теперь $z=(\beta,\alpha,-1)$, где $\beta\in\mathbb{C}$ --- любое комплексное число.
Поскольку $|\alpha|=1$, отсюда следует, что
$$
2-2\Re \{{\bar{\alpha f(x-y)}}\} \geq -|\beta|^2 -
2 \Re \{{[\alpha f(x)-f(y)]\bar{{\beta}}}\}.
$$
Положим теперь $\beta=-[\alpha f(x)-f(y)]$. Тогда
$$
2(1-\Re \{{\alpha f(x-y)})\} \ge |\alpha f(x)-f(y)|^2.
$$

\begin{remark}
Доказательство теоремы проведено по той же схеме, по которой М.Г. Крейн доказывал неравенство \eqref{eq2} в 1943 году. Неравенство \eqref{eq3} также следует из формально более общего неравенства (iv) из книги Сасвари \cite{Zas}, теорема 1.4.12.
\end{remark}

При $\alpha=-1$ из \eqref{eq3} получаем неравенство для п.о.ф.
\begin{equation}\label{eq4}
|f(x)+f(y)|^2 \le 2 f(0) [f(0)+\Re f(x-y)].
\end{equation}
Отсюда следует, что если в некоторой точке $T\neq 0$ справедливо равенство $f(T)=-f(0)$, то
$$
f(x+T)=-f(x), f(x+2T)=f(x), x\in\mathbb{R}.
$$
Аналогичное комплексное следствие вытекает из теоремы 1.

\begin{corollary}
Если в некоторой точке $T\neq 0$ выполнено равенство $f(T)=\alpha f(0)$, где $|\alpha=1$, то справедливо тождество
\begin{equation}\label{eq5}
f(x+T)=\alpha f(x), x\in\mathbb{R}.
\end{equation}
\end{corollary}

Доказательство. Возьмём в \eqref{eq3} $y=x+T$. Имеем
\begin{equation}\label{eq6}
|\alpha f(x) - f(x+T)|^2 \le 2 f(0) \Re[f(0)-  {\alpha f(-T)}].
\end{equation}
Но $\alpha f(-T)=\alpha \bar{f(T)}=\alpha\bar{\alpha}f(0)=f(0)$, значит, правая часть \eqref{eq6} равна нулю, и, следовательно, выполнено \eqref{eq5}.

Следствие 1 можно найти в (\cite{Gor1}, лемма 1), где используется теорема Бохнера. У нас это следствие получается прямо из неравенства \eqref{eq3}.

Классическое неравенство \eqref{eq2} получено в статье М.Г. Крейна \cite{Kre} 1943 года. В этой статье рассматривается проблема продолжения функции, положительно определённой на $(-R,R)$, на всю вещественную ось.

Возникает вопрос: если две п.о.ф. совпадают на некотором интервале $(-R,R)$, то будут ли они совпадать на всей вещественной оси? Ответ отрицательный, как показывает следующий

Пример. Рассмотрим функцию --- "домик" (tent function)

$$
f(x)=\begin{cases}
2-|x| & \text{ при $|x|\le 2$}, \\
0     & \text{ при $|x|> 2$}.
\end{cases}
$$
Известно, что она положительно определена (см. \cite{Bha}, p. 149). Рассмотрим также "маленький домик"
$$
g_1(x)=\begin{cases}
1-|x| & \text{ при $|x|\le 1$}, \\
0     & \text{ при $|x|> 1$}.
\end{cases}
$$
и постоянную функцию $g_2(x)=1$. Тогда $g(x)=g_1(x)+g_2(x)$ будет п.о.ф., и $f(x)=g(x)$ при $x\in[-1;1]$. Но на всей оси $f$ и $g$ не совпадают.

\section{Неравенство Ю.В. Линника и его усиления }

Это неравенство справедливо для вещественных п.о.ф., которые в этом параграфе будем обозначать $u(x)$. Такие функции совпадают с известными в теории вероятностей характеристическими функциями симметричных распределений \cite{Lin}--\cite{Fel}. В книге \cite{Lin} установлена

\begin{theorem}
Справедливо неравенство

\begin{equation}\label{eq7}
u(0) - u(2x) \le 4(u(0) - u(x), x\in\mathbb{R}.
\end{equation}
\end{theorem}

Доказательство. Известно \cite{Bha}, что $|u(x)|\le u(0)$. Если $u(0)=0$, то \eqref{eq7} тривиально. Если же $u(0)>0$, то можно считать, что $u(0)=1$. Положим в \eqref{eq4} $y=-x$ и учтём, что $u(-x)=u(x)$. Получим $4 u^2(x)\le 2(1+u(2x)$, или после очевидных преобразований
$$
1+ 2u^2(x) \le 2+u(2x),
$$
$$
1-u(2x) \le 2(1-u^2(x)) \le 4(1-u(x)).
$$
Последнее неравенство равносильно $(u-1)^2 \ge 0$.

Неравенство \eqref{eq7} приведено также в книге В. Феллера (\cite{Fel}, с. 560) без упоминания Ю.В. Линника. В отличие от \cite{Lin}--\cite{Fel} в нашем доказательстве не используется теорема С. Бохнера об интегральном представлении п.о.ф. Попутно получено неравенство
\begin{equation}\label{eq8}
1-u(2x) \le 2(1-u^2(x)),
\end{equation}
которое обращается в равенство для п.о.ф. $u(x)=\cos(x)$.

Константа $4$ в неравенстве \eqref{eq7} является наилучшей. Действительно, возьмём функцию $u(x)=\exp(-x^2)$, запишем для неё неравенство \eqref{eq7}, поделим на $x^2$ и перейдём к пределу при $x\to 0$. Получим неравенство $4\le 4$.

Теперь рассмотрим следствия неравенства Ю.В.Линника в случае $u(0)=1$.

В этом случае $|u(x)| \le 1$ для всех $x$ и неравенство \eqref{eq7} принимает вид
\begin{equation} \label{eq9}
1 - u(2x) \le 4 (1 - u(x)).
\end{equation}
Элементарно проверяется, что неравенство \eqref{eq9} равносильно неравенству
\begin{equation} \label{eq10}
1 + u(x) \le \frac{7 + u(2x)}{4}.
\end{equation}
Многократно повторяя неравенство \eqref{eq9}, получим неравенство
\begin{equation} \label{eq11}
1 - u(2^m x) \le 4(1 - u(2^{m-1} x)) \le 4^m (1 - u(x)), m\in \mathbb{N}.
\end{equation}
Наконец, если использовать неравенства \eqref{eq8} и \eqref{eq10}, то выражение $1 - u(2^m x)$ можно оценить более точно. Покажем это, например, для $m=3$:

$$
1 - u(8x) \le 2(1 - u^2(4x))=2(1-u(4x))(1+u(4x))\le
$$
$$
\le 4(1 - u^2(2x)) \frac{7+u(8x)}{4} = 4(1-u(2x))(1+u(2x))\cdot\frac{7+u(8x)}{4} \le
$$
$$
\le 8(1 - u^2(2x))\cdot \frac{7+u(4x)}{4}\cdot \frac{7+u(8x)}{4} \le
$$
$$
\le 8(1 - u(x)) \cdot \frac{7+u(2x)}{4} \cdot \frac{7+u(4x)}{4} \cdot \frac{7+u(8x)}{4}.
$$

В общем случае получаем такое уточнение неравенства Ю.В.Линника
\begin{equation} \label{Lin}
1 - u(2^m x) \le 2^m (1 - u(x)) \prod_{k=1}^m \frac{7+u(2^k x)}{4}, m\in\mathbb{N}.
\end{equation}
Это неравенство лучше, чем \eqref{eq11}, так как каждый сомножитель в произведении не превосходит числа 2.

\section{Основные многоточечные неравенства}

Для вещественных п.о.ф. $u(x)$ можно установить неравенства, в которых участвуют суммы произвольных значений $x_1, \dots, x_n \in\mathbb{R}.$

\begin{theorem}
Для вещественных непрерывных п.о.ф. $u(x)$ справедливо неравенство
\begin{equation}\label{eq12}
u(0) - u(x_1+\dots+x_n) \le n \sum_{k=1}^n (u(0)-u(x_k))
\end{equation}
для всех $n$ и любых $x_1, \dots, x_n \in\mathbb{R}.$
\end{theorem}

Доказательство. Можно считать, что $u(0)=1$. По теореме Бохнера неравенство \eqref{eq12} переписывается в виде
\begin{equation}\label{eq13}
\int_{-\infty}^{+infty} \left[1 - \cos(t(x_1+\dots+x_n))\right]\,\mu(dt) \le n \sum_{k=1}^n
\int_{-\infty}^{+infty} (1 - \cos(t x_k))\,\mu(dt),
\end{equation}
где $\mu$ --- некоторая вероятностная мера. Выполнение \eqref{eq13} для всех допустимых мер $\mu$ равносильно выполнению неравенства
\begin{equation}\label{eq14}
1 - \cos(t(x_1+\dots+x_n)) \le n \sum_{k=1}^n (1 - \cos(t x_k))
\end{equation}
для всех $t, x_1, \dots, x_n \in\mathbb{R}$. Положим $s_k=tx_k/2$ и преобразуем \eqref{eq14} к виду
\begin{equation}\label{eq15}
\sin^2(s_1+\dots+s_n) \le n \sum_{k=1}^n \sin^2(s_k).
\end{equation}
Последнее неравенство следует из другого тригонометрического неравенства
\begin{equation}\label{eq16}
|\sin(s_1+\dots+s_n)| \le \sum_{k=1}^n |\sin(s_k)|,
\end{equation}
которое легко доказывается индукцией по $n$, после чего для доказательства неравенства \eqref{eq15} остаётся возвести \eqref{eq16} в квадрат и применить неравенство Коши--Буняковского. Теорема доказана.

\begin{corollary}
 Пусть $n$ --- нечётное число, $f$ --- непрерывная комплексная п.о.ф. Тогда для любых $x_1,\dots, x_n, y_1, \dots, y_n \in\mathbb{R}$ справедливо неравенство
\begin{equation}\label{eq17}
|f(x_1+\dots+x_n) - f(y_1+\dots+y_n)|^2 \le 2n f(0) \sum_{k=1}^n [f(0)-\Re f(x_k - y_k)].
\end{equation}
\end{corollary}
\begin{remark} Это неравенство получено в статье Е.А.Горина (\cite{Gor1}, теорема 1).
\end{remark}

Доказательство следствия 2. Считаем, что $f(0)=1$. По неравенству Крейна \eqref{eq2} левая часть $L$ неравенства \eqref{eq17} оценивается так:
$$
L \le 2(1 - \Re f(\sum_{k=1}^n (x_k - y_k)).
$$
применим неравенство \eqref{eq12} к п.о.ф. $u(x)= \Re f(x)$. Получим
$$
L \le 2n \sum_{k=1}^n (1 - u(x_k - y_k)),
$$
что и требовалось доказать.

Как и в  неравенстве Крейна в основном неравенстве \eqref{eq12} знак "минус"\, можно поменять на знак "плюс". Тогда получим два новых неравенства.

\begin{theorem} Для непрерывной вещественной п.о.ф. $u(x)$ справедливо неравенство
\begin{equation} \label{eq18}
u(0) - u(x_1+\dots+x_n) \le n \sum_{k=1}^n (1 + u(x_k))
\end{equation}
для любого чётного $n$ и $x_1,\dots,x_n \in\mathbb{R}$.
\end{theorem}

Доказательство. Совершенно так же, как в теореме 3 неравенство \eqref{eq18} сводится к тригонометрическому неравенству
$$
\sin^2(s_1+\dots+s_n) \le n \sum_{k=1}^n \cos^2{s_k}.
$$
Для этого надо положить $s_k=\pi/2 - t_k$ и учесть чётность $n$. Тогда требуемое неравенство сведётся к уже доказанному неравенству \eqref{eq15}.

Отметим, что при нечётном $n$ неравенство  \eqref{eq18} не выполняется, например, для функции $u(x)=\cos(x)$ и при выборе точек $x_k=\pi$.

Для комплексной п.о.ф. $f(x)$ из  \eqref{eq18} получаем неравенство
\begin{equation}\label{eq19}
|f(x_1+\dots+x_n) - f(y_1+\dots+y_n)|^2 \le 2n f(0) \sum_{k=1}^n [f(0)+\Re f(x_k - y_k)],
\end{equation}
верное при любом чётном $n$.

Это неравенство следует из основной леммы Е.А. Горина \cite{Gor2}  при  $\zeta=-1$  и четном  $n$.

Отметим, что неравенство с разными знаками "плюс"/"минус"\  могут нарушаться при всех $n$. А вот неравенство с двумя плюсами выполняется при нечётном $n$
\begin{equation}\label{eq20}
u(0) + u(x_1+\dots+x_n) \le n  \sum_{k=1}^n [1+u(x_k)]
\end{equation}
при всех $x_1,\dots,x_n \in\mathbb{R}$. Действительно, оно сводится ещё к одному  тригонометрическому неравенству
$$
\cos^2(s_1+\dots+s_n) \le n \sum_{k=1}^n \cos^2(s_k),
$$
которое заменой $s_k=\frac{Pi}{2} - r_k$ при нечётных $n$ переходит в \eqref{eq15}.

Для комплексной п.о.ф. f(x) справедливо неравенство
\begin{equation}\label{eq21}
|f(x_1+\dots+x_n) + f(y_1+\dots+y_n)|^2 \le 2n f(0) \sum_{k=1}^n [f(0)+\Re f(x_k - y_k)],
\end{equation}
для любого нечётного $n$. Для доказательства надо применить сначала неравенство \eqref{eq4}, а затем \eqref{eq20}. При чётном $n$ неравенства \eqref{eq20}--\eqref{eq21} не выполняются, что видно из примера $u(x)=f(x)=\cos(x), x_k=\pi, y_k=0$.

В заключение отметим, что для полученных в работе неравенств возможно их дальнейшее уточнение с использованием методики обобщений неравенств Коши--Буняковского, развитой в \cite{Sit1}--\cite{Sit2}. Полученные результаты могут быть также применены к получению неравенств для специальных функций \cite{KaSi1}--\cite{KaSi2}, некоторых интегральных неравенств \cite{PeKo}, в теории сплайнов \cite{IgPe} и для обоснования однозначной разрешимости задач аппроксимации функций квадратичными экспонентами \cite{Sit3}--\cite{Sit4}.

\renewcommand\refname{\bf Литература}

\bibliographystyle{amsplain}

\vskip 20mm

АВТОРЫ:\\
\\
\\
Ситник С.М.\\
Воронежский институт МВД, Воронеж, Россия.\\
Sitnik S.M., Voronezh Institute of the Ministry of Internal Affairs of Russia.\\
mathsms@yandex.ru, pochtasms@gmail.com
\\
\\
Певный А.Б.\\
Сыктывкарский госуниверситет, Сыктывкар, Россия.\\
Pevnyi A.B., Syktyvkar State University, Syktyvkar, Russia.\\
pevnyi@syktsu.ru
\end{document}